\title{\vspace{-0.7cm} Local resilience of graphs }

\author{
Benny Sudakov \thanks{
Department of Mathematics, UCLA, Los Angeles, CA 90095 and
Institute for Advanced Study, Princeton. E-mail:
bsudakov@math.princeton.edu.
Research supported in part by NSF CAREER award DMS-0546523, NSF grant
DMS-0355497, USA-Israeli BSF grant, Alfred P. Sloan fellowship, and
the State of New Jersey.}
\and V. H. Vu
\thanks{
Department of Mathematics, Rutgers University, Piscataway, NJ 08854, USA.
E-mail: vanvu@math.rutgers.edu. Research supported in part by
NSF CAREER award.
}
    }

\documentclass [11pt]{article}
\usepackage{amsfonts}
\usepackage{amsmath}
\usepackage{latexsym}
\newtheorem{theo}{Theorem}[section]
\newtheorem{theorem}[theo]{Theorem}

\newtheorem{lemma}[theo]{Lemma}

\newtheorem{fact}[theo]{Fact}
\newtheorem{defn}[theo]{Definition}

\newtheorem{problem}[theo]{Problem}

\def\CP{{\mathcal P}}
\def\ep{{\epsilon}}

\date{}

\begin{document}
\maketitle

\begin{abstract}
In this paper, we initiate  a systematic study  of graph
resilience. The (local) resilience  of a graph $G$ with respect to
a property $\CP$ measures how much one has to change $G$ (locally)
in order to destroy $\CP$. Estimating the resilience  leads to many new and challenging problems.
Here we focus on random and pseudo-random graphs and prove several sharp results.
\end{abstract}
\section {The notion of resilience }

A typical result in graph theory is of the following form:

\vskip2mm

 \centerline {\it A  graph  $G$ (from certain class)  possesses a property $\CP$. }

\vskip2mm \noindent In this paper, we would like to  investigate
the following general problem:

\vskip2mm

\centerline {\it How strongly does  $G$ posses $\CP$?}
\vskip2mm

To study  this question, we define {\it the
resilience of $G$ with respect to $\CP$}, which measures how much
one should change $G$ in order to destroy $\CP$.
There are two natural kinds of resilience: global and local. It is
more convenient to first define these quantities with respect to
monotone increasing properties ($\CP$ is monotone increasing if it
is preserved under edge addition).

\begin{defn} (Global resilience)
Let $\CP$ be an increasing monotone property. The global
resilience  of $G$ with respect to $\CP$ is the minimum number $r$
such that by deleting $r$ edges from $G$ one can obtain a graph not having
$\CP$.
\end{defn}

\noindent The notion of global resilience  is not new. In fact,
problems about global resilience  are  popular in extremal graph
theory. For example, the celebrated theorem of Tur\'an \cite{Tu}
gives the answer to the following question:

\vskip2mm

\centerline{\it How many edges should one delete from the complete graph $K_n$ to make it $K_k$-free?
}

\vskip2mm

The main focus of this paper is on the {\it local resilience },
which eventually leads to a host of intriguing  new questions. To
start, let us notice that one can destroy many properties by
simple local changes.  For instance, to destroy the hamiltonicity,
it suffices to delete all edges adjacent to one vertex. This
motivates the following notion.

\vskip2mm

\begin{defn} (Local resilience)
Given a monotone increasing property $\CP$. The local resilience
of a graph $G$ with respect to $\CP$ is the minimum number $r$ such
that by deleting at each vertex of $G$ at
most $r$ edges one can obtain a graph not having $\CP$.
\end{defn}

\noindent 
Note that a classical theorem of Dirac (see, e.g.,
\cite{D}) can be formulated in this language as follows:

\vskip2mm \centerline{\em The local resilience of $K_n$ with
respect to hamiltonicity is $ \lfloor n/2 \rfloor$.} \vskip2mm

\noindent If $\CP$ is not monotone (for instance, having a
non-trivial automorphism is such a property), then we may have to
both delete and add edges. This leads  to the general definition.

\begin{defn}
Given a property $\CP$. The local resilience  of $G$ with respect
to $\CP$ is the minimum number $r$ so that there is a graph $H$ on
$V(G)$ with maximum degree at most $r$ such that the graph $G
\triangle H$ does not have $\CP$. \end{defn}

One can observe that there is a certain duality between properties
and resilience. If the property is local (such as containing a
triangle), then it makes more sense to talk about the global
resilience. If the property is global (such as being
hamiltonian), then the local resilience  seems to be the right
parameter to consider.

\vskip2mm

The notion of local resilience was first introduced few years ago
in a paper of Kim and Vu \cite{KV} for a purpose which we will
discuss in concluding remarks. It was also used implicitly by Kim,
Sudakov and Vu in \cite{KSV}  to prove a conjecture of Wormald
(see the subsection on Symmetry in the next section). We are now
convinced that this notion is of independent interest, and, with
this paper, aim to initiate a systematic investigation.

\vskip2mm
To generate new questions, one can match any
interesting graph with any natural property and ask for the
corresponding resilience. Here is a random example:

\vskip2mm

{\it What is the local resilience  of the hypercube with respect
to containing a perfect matching?}

\vskip2mm

\noindent In other words,  what is the minimum degree of a subgraph of  the hypercube
 containing  an edge from every  perfect matching?

 \vskip2mm

 The notion of resilience  is not restricted to graphs. In fact,
one can define it for virtually any structure. For instance, one
can consider the resilience of matrices. The resilience, in this
case,  tells us how many entries need to be change in order to
destroy a certain property. Here is an example, suggested by Kahn
(private communication):

\vskip2mm

{\it What is the global resilience  of the  Hadamard matrix with
respect to being non-singular?}

\vskip2mm

\noindent In other words, how many entries of a Hadamard matrix
needs to be changed in order to make the matrix singular? The
problem seems even more interesting if we restrict the changes to
a small set of values (e.g., a new entry can only take values in
$\{-1,1\}$ or $\{0,-1,1\}$).

\vskip2mm

Instead of producing a long list of  problems here (we leave this
task to the interested readers), we would like to mention that the
notion of resilience not only generates new questions, but also
serves as a useful tool for proving new results. We will discuss
an example in Section  2.4 and another in Section 6.1. As a matter
of fact, these applications were the original motivation for
developing this notion.

\vskip2mm

To conclude this section, let us point out that  the notion of
resilience may have a relation to the theory of communication. A
basic scenario in this theory involves three parties. A sender
(Alice), a receiver (Bob) and an eavesdropper (Charles). Alice
wants to send a message to Bob. On the other hand, Charles
(usually  an adversary) can (and want to) intercept and change the
message.

Now imagine that the message has a form of a binary matrix, which is the adjacency matrix
of a graph. The information Bob wants to learn from the message is hidden in a property of
the graph. Can Charles destroy this information?

The answer depends on how much power Charles possesses. If his
power was unlimited, then he could simply turn all bits to zero
and perhaps destroy any interesting information the message may
contain. But what if his power is limited? In fact, it is a
common assumption that the matrix is  sent row by row and Charles
could only change a few bits in each row. It is now clear that if
the original graph has high local resilience (with respect to the
particular property Bob is interested in), then Charles cannot do
anything to prevent Bob from learning the truth.

\section{The resilience  of random graphs}

 Fix a property $\CP$ and consider the set of all (labeled) graphs on the same large vertex set.
The local resilience is viewed as  a graph function $f_{\CP}(G)$.
The value of this function can change significantly from graph to
graph, and for some graphs could be very hard to determine. Thus,
a natural thing to do is to study the typical behavior of this
function in a large set of graphs. This, as one may expect, leads
immediately to the consideration of the local resilience  of
random graph.

The most commonly used model of random graphs, sometimes synonymous with the term
``random graph", is the so called binomial random graph $G(n,p)$ introduced by
Erd\H os and R\'enyi. The random graph $G(n,p)$ denotes the probability space whose
points are graphs on a fixed vertex set $[n]=\{1,\dots, n\}$
where each pair of vertices forms an edge, randomly and independently, with
probability $p$. Here $p$ is a positive number at most one and can depend on $n$.
The above operation generates a distribution on the set of all (labeled)
graphs on $n$ vertices. Abusing the notation slightly, we call this distribution the
$G(n,p)$ distribution. It is well known, and easy to prove, that this distribution
concentrates on the graphs with roughly ${n \choose 2} p$ edges.
We say that the random graph $G(n,p)$
possesses a graph property $\CP$ {\em almost surely},
or a.s. for short, if the probability that $G(n,p)$ satisfies $\cal P$
tends to 1 as the number of vertices $n$ tends to infinity.

\begin{defn}
Consider random graph $G(n,p)$ and a fixed property $\CP$. The
local resilience  of $G(n,p)$ with respect to $\CP$ is the minimum
number $r$ so that almost surely there is a graph $H$ on $[n]$
with maximum degree at most $r$ such that the graph $G(n,p)
\triangle H$ does not have $\CP$.
\end{defn}

\vskip2mm \noindent {\bf Remark.}\, It is important to notice that
after one has sampled a graph from the distribution $G(n,p)$, the
adversary is allowed to find graph $H$ with the smallest possible 
maximum degree in order to destroy the property $\CP$. The graph $H$ can change from
sample to sample.

\vskip2.5mm \noindent {\bf Remark.}\, Instead of considering the
distribution $G(n,p)$, one can consider the set of all graphs on
$n$ vertices with ${n \choose 2} p$ edges, equipped with the
uniform distribution. All results we prove for $G(n,p)$ will hold
with respect to this distribution as well.

\vskip2.5mm \noindent 
Next, we will describe several estimates which we obtained for the
local resilience of random and pseudo-random graphs with respect
to various graph properties. Some of these results are
asymptotically best possible. Throughout this paper, $\log$ denotes logarithm in the natural base $e$.
For two functions $f(n)$ and $g(n)$ the notation $f \gg g$ means that
$f/g \rightarrow \infty$ together with $n$.

\subsection {Perfect matching}

Let $G$ be a graph on $n$ vertices, where $n$ is even.
A perfect matching of $G$ is a set of disjoint edges which together cover all $n$
vertices. A classical theorem of Erd\H{o}s and R\'enyi \cite{ER2} shows that for any fixed $\ep >0$
if $p > (1+\ep) \log n/n$, then almost surely $G(n,p)$ has a perfect
matching.

Let us now give a lower bound for the local resilience  of the
property of having a perfect matching. One natural way to destroy
all perfect matchings is the following. Split the vertex set of
$G$ into two parts $X$ and $Y$ of size $n/2+1$ and $n/2-1$
respectively. Then delete all edges inside the set $X$. Thus $X$
becomes an independent set  and it is impossible to match all of
its vertices  with vertices of $Y$, since $|Y|<|X|$. In other
words, if we define $f_1(G)$ to be the minimum of the maximum
degree of the subgraph of $G$ induced by a subset of size $n/2+1$
then

\begin{fact} The local resilience  of $G$ with respect to having a perfect matching is
at most
$f_1(G)$.
\end{fact}

In $G(n,p)$ (with $p$ sufficiently large), with high probability
all vertices have degree $(1+o(1))np$. Thus,  one would expect
that almost surely random graph has an induced subgraph
on $n/2+1$ vertices whose maximum degree is $(1/2+o(1))np$. 
So, the local resilience  of
$G(n,p)$ with respect to having a perfect matching is at most
$(1/2+o(1))np$. Our first theorem shows that this trivial upper
bound is actually the truth.

\begin{theorem}
\label{theorem:Perfect1} For even $n$ and $p \gg \log n/n$, the
local resilience of $G(n,p)$ with respect to having a perfect
matching is $(1/2+o(1))np$.
\end{theorem}

Roughly speaking, this theorem tells us that for most graphs on $n$ vertices with
approximately
${n \choose 2}p$ edges, there is no essentially better way to remove all perfect
matchings than the one considered above. We will present a more precise statement of
Theorem \ref{theorem:Perfect1} together with its proof in Section 3.

For $n$ odd, one can talk about  nearly perfect matchings (which
leave one vertex unmatched). Theorem \ref{theorem:Perfect1} also
holds in this case with essentially the same proof.

\subsection{Hamiltonicity}

A graph $G$ is hamiltonian if it contains a cycle going through
all vertices. Another classical theorem in random graph theory
which was proved by Koml\'os and Szemer\'edi \cite{KomSze83} and
by Korshunov \cite{Kor} (extending earlier result of P\'osa
\cite{Pos76}) shows that if $p > (1+\ep) \log n/n$, then almost
surely $G(n,p)$ is hamiltonian.

Following the ideas in the previous section, we can remove all
hamiltonian cycles of $G$ by splitting the vertex set of $G$ into
two parts whose sizes differ by at most two and deleting all the
edges inside the larger part. Then, similarly to the previous
section we have that

\begin{fact} The local resilience  of $G$ with respect to having a Hamilton cycle is at
most $f_1(G)$. \end{fact}

Another construction is to split the vertex set into two parts
whose sizes differ by at most one  and delete all edges between
them. This, however, does not change the asymptotic of the lower
bound with respect to a random graph.
Therefore for sufficiently large $p$, we can again
conclude  that the local resilience  of $G(n,p)$ with respect to
being hamiltonian is at most $(1/2+o(1))np$. Here we prove the
matching lower bound.

\begin{theorem} \label{theorem:Hal1} For all $p > \log ^{4} n/n$, the local
resilience  of $G(n,p)$ with respect to being hamiltonian is
$(1/2+o(1))np$.
\end{theorem}

\noindent
A more accurate version of this theorem and its proof will be
presented in Section 4. Our proof relies only on properties of the
edge distribution of $G(n,p)$ and therefore can be extended to
pseudo-random graphs.

Consider a graph $G$. Let  $\lambda_1(G) \geq \lambda_2 (G) \geq
\ldots \geq \lambda_n (G)$ be the eigenvalues of its adjacency
matrix. The quantity $\lambda(G)= \max_{i \geq 2} |\lambda_i(G)|$
is called the {\em second eigenvalue} of $G$ and plays an
important role. We say that $G$  is  an {\em
$(n,D,\lambda)$-graph} if it is $D$-regular, has $n$ vertices and
$\lambda(G)$ is at most $\lambda$. It is well known (see, e.g.,
\cite{AloSpe00} for more details) that if $\lambda$ is much
smaller than the degree $D$, then $G$ has strong pseudo-random
properties, i.e., the edges of $G$ are distributed like in random
graph $G(n,D/n)$. In Section 4 we will show that if
$D/\lambda>\log^2 n$ then the local resilience  of any
$(n,D,\lambda)$-graph with respect to being hamiltonian is
$(1/2+o(1))D$.

\subsection{Chromatic number}

A graph $G$ is {\em $k$-colorable} if one can assign $k$ different
colors to the vertices of $G$ so that the colors of adjacent
vertices are different. The {\em chromatic number} of $G$ is the
smallest $k$ such that $G$ is $k$-colorable. We use the standard
notation $\chi (G)$ for the chromatic number of $G$. It was proved
by Bollob\'as \cite{Bo} for large $p$ and by \L uczak \cite{Lu}
for small $p$  that almost surely
$\chi\big(G(n,p)\big)=(1+o(1))\frac{n}{2\log_b (np)}$, where
$b=1/(1-p)$.

In Section 5, we consider the local resilience  of $G(n,p)$ with
respect to the colorability property. The task is to determine the
largest $r$ such that for almost all samples $G$ from $G(n,p)$,
the following holds. By adding at most $r$ edges to each vertex of
$G$, one cannot increase the chromatic number of $G(n,p)$ by more
than $o(\chi(G(n,p))$. Trivially $r \le (1+\ep) \chi(G(n,p))$ for
any fixed $\ep > 0$, since we can add a clique of that size and
the chromatic number of the graph will grow by factor at least
$1+\ep$. From below, we have the following bound.

\begin{theorem} Let $p$ be a fixed, small positive number.
Assume that $n^{-1/3+\ep} \le p \le 3/4$. The local resilience  of
$G(n,p)$ with respect to being $(1+o(1)) \frac{n}{2\log_b
(np)}$-colorable is at least $np^2/ \log^5 n$.
\end{theorem}

\noindent Note that for the uniform case when $p=1/2$, this lower
bound is off by only  polylogarithmic factor. For $p$ below
$n^{-1/3 +\ep}$, we have the following weaker result.

\begin{theorem}
\label{t1}
For every positive integer $d$ and for every $\ep >0$ there is a constant
$c=c(d,\epsilon) $ such that the following hold. For any $p > c/n$, almost surely

$$\max_{H, \Delta (H) \le d} \chi ((G(n,p) \cup H) \le (1+\ep) \chi (G(n,p)).
$$\end{theorem}

Informally, this result says that for any given $d$ and for most
graphs on $n$ vertices with average degree at least $c$ (which is
sufficiently large compared to $d$), adding any graph of maximum
degree $d$ has very little impact on the chromatic number. It is
best to compare this theorem against the following folklore result (see,
e.g., \cite{Lovaszbook} Chapter 9).

\begin{fact} Let $G$ and $H$ be two graphs on the same set of points. Then
$$\chi (G \cup H) \le \chi (G) \chi (H),$$
and there are pairs of $G$ and $H$ such that the equality  holds.
\end{fact}

\subsection {Symmetry}

Another classical theorem of Erd\H os and R\'enyi \cite{ER1}
states that for $p \ge (1+\ep) \log n /n $,  $G(n,p)$ is almost
surely non-symmetric, i.e., has no non-trivial automorphisms. The
question here is how much should one change the graph to make it
symmetric. The global resilience  with respect to being
non-symmetric was determined in \cite{ER1}.

\begin{theorem} \label{theorem:ER} If both $1-p, p \gg \log n/n$, then the global
resilience  of $G(n,p)$  with respect to being non-symmetric is
$(2+o(1))np(1-p) $.
\end{theorem}

To explain the quantity $(2+o(1))np(1-p) $, notice that this is
the number of edges one needs to delete to obtain two vertices
with the same neighborhoods. The transposition of such two
vertices is a non-trivial automorphism of the (modified) graph.
Notice that in this case, one needs to delete $(1+o(1))np(1-p)$
edges from each vertex. In \cite{KSV}, the authors, together with
 Kim, proved a strengthening of Theorem \ref{theorem:ER} which
implies the following result.

\begin{theorem} If both $1-p, p \gg \log n/n$, then the local resilience
of $G(n,p)$ with respect to being non-symmetric
is $(1+o(1))np(1-p) $.
\end{theorem}

As a corollary, we proved in \cite{KSV} that random regular graphs
(with large degrees) are almost surely non-symmetric, confirming a
conjecture of Wormald. The interesting point here is that local
resilience  was not the main object under study in \cite{KSV}, but
it was used as a tool  to prove a new result. The main idea in
\cite{KSV} is as follows. Consider the indicator graph function
$I(G)$, where $I(G)=1$ if $G$ is non-symmetric and $0 $ otherwise. 
We want to show that with high probability $I(G)=1$ where
$G$ is a random regular graph. One may want to view this
statement as a sharp concentration result, namely, $I(G)$ is a.s.
close to its mean. However, it is impossible to
prove a sharp concentration result for a  random variable having
only two values close to each other. The idea here is to ``blow up"
$I(G)$ using the notion of local resilience. Instead of $I(G)$ we
used a function $D(G)$ which (roughly speaking) equals the local
resilience  of $G$ with respect to being non-symmetric. This
function is zero if $G$ is symmetric and rather large otherwise.
This gives us room to show that $D(G)$ is strongly concentrated around a large
positive value, and from this we can conclude that almost surely
the random regular graph is non-symmetric.

\vspace{0.3cm}
\noindent
{\bf Notation.}\,  We usually write $G=(V,E)$ for a graph $G$ with vertex
set $V=V(G)$ and edge set $E=E(G)$, setting $n=|V|$ and
$e(G)=|E(G)|$. If $X \subset V$ is a subset of the vertex set then
$e_G(X)$ denotes the number of those edges of $G$ with both
endpoints in $X$. Similarly, we write $e_G(X,Y)$ for the number of
edges with one endpoint in $X$ and the other in $Y$. We denote by
$N_G(v)$ the set of vertices adjacent to a vertex $v$ and by
$d_G(v)=|N(v)|$ the degree of $v$. More generally, for a subset $X
\subset V$, $N_G(X)$ is the set of all neighbors of vertices of
$X$ in $G$. Note that $N_G(X)$ is not necessarily disjoint
from $X$.

To simplify the presentation, we often omit
floor and ceiling signs whenever these are not crucial and make no attempts to optimize the
absolute constants. 

\section{Perfect matchings}

The following is a more precise version of Theorem
\ref{theorem:Perfect1}. To deduce Theorem \ref{theorem:Perfect1}
from it, simply notice that if $p \gg \log n /n$, then the error
term $8 \sqrt{np \log n } = o(np) $.

\begin{theo}
\label{matching}
For even $n$ and $p \geq 400\log n/n$ the random graph $G(n,p)$
almost surely has the following properties:

\noindent
$(i)$\, If  $H$ is a subgraph of $G=G(n,p)$ with maximum degree
$$\Delta(H) \leq np/2-8\sqrt{np \log n},$$
then $G'=G-H$ has a perfect matching.

\noindent
$(ii)$\, There exist a subgraph $H$ of  $G=G(n,p)$ with maximum degree
$$\Delta(H) \leq np/2+2\sqrt{np \log n}$$
such that $G'=G-H$ has no perfect matching.
\end{theo}

To prove the theorem we need the following simple lemma about properties of
random graphs.

\begin{lemma}
\label{properties}
If $p \geq 400\log n/n$ then almost surely,

\noindent
$(i)$\, The minimum degree of $G(n,p)$ is at least $np-2\sqrt{np \log n}$,

\noindent
$(ii)$\, The maximum degree of $G(n/2+1,p)$ is at most $np/2+2\sqrt{np \log n}$,

\noindent
(iii)\, The number of edges between any two disjoint subsets $A, B$ of $G(n,p)$ of size
$|A|=|B|=s \leq n/4$ is at most
$$e(A,B) \leq s \Big( \frac{np}{4}+\sqrt{2np \log n}\, \Big).$$
\end{lemma}

\noindent {\bf Proof.}\, (i)\, Since the degree of every vertex $v$ of $G(n,p)$ is
binomially distributed random variable with parameters $n-1$ and $p$, it follows by the
standard Chernoff-type estimates (see, e.g., Chapter 2.1 in \cite{JanLucRuc00}) that
$$\mathbb{P}\Big[d(v)-(n-1)p \leq -\sqrt{3np \log n}=-t\Big] \leq
e^{-\frac{t^2}{2np}}=n^{-3/2}.$$
Therefore the probability that $G(n,p)$ has a vertex of
degree less than $np-2\sqrt{np \log n}<(n-1)p-\sqrt{3np \log n}$ is at most $n^{-1/2}=o(1)$.

\vspace{0.1cm}
\noindent
(ii)\,  Let $v$ be a vertex of $G(n/2+1,p)$ and let $t=2\sqrt{np \log n}$.
Since $2\sqrt{np \log n}< np$, by Chernoff-type estimate we have that
$$\mathbb{P}\Big[d(v)-np/2\geq t\Big] \leq e^{-\frac{t^2}{2(np/2+t/3)}} \leq
e^{-\frac{t^2}{2np}}= n^{-2}.$$
Thus the probability that $G(n/2+1,p)$ has a vertex of degree larger than
$np/2+2\sqrt{np \log n}$ is $o(1)$.

\vspace{0.1cm}
\noindent
(iii) The number of edges between two disjoint sets $A, B$ of size $s$ is binomially distributed and has
expectation $\lambda=s^2p$. Set $t=s\sqrt{2np \log n}$. Since $s\leq n/4$, by Chernoff-type
estimates, we have that
$$\mathbb{P}\Big[e(A,B) \geq s \big(np/4+\sqrt{2np \log n}\,\big)\Big] \leq
\mathbb{P}\Big[e(A,B) \geq s^2p+s\sqrt{2np \log n}=\lambda +t\Big] \leq e^{-\frac{t^2}{2(\lambda
+t/3)}}.$$
As $p \geq 400\log n/n$ and $s \leq n/4$, we have that $\sqrt{2np \log n}<np/10$
and $2sp\leq np/2$. Therefore
$$\frac{t^2}{2(\lambda +t/3)}>s \cdot \frac{2np \log n}{2sp+\sqrt{2np \log n}}\geq
s \cdot \frac{2np \log n}{np/2+np/10}>3s\log n.$$
This implies that the probability that  there are two disjoint sets which violates
the assertion of the lemma is at most
$$ \sum_{s\leq n/4} {n \choose s} {n \choose s} e^{-3 s\log n} \leq
\sum_{s\leq n/4} n^{2s}\, n^{-3s}=\sum_{s\leq n/4} n^{-s}=o(1).$$
This completes the proof.
\hfill $\Box$

\vspace{0.2cm}
\noindent
{\bf Proof of Theorem \ref{matching}.}\, (i)\, Let $G$ be a graph
of order $n$ with minimum degree at least $np-2\sqrt{np \log n}$ such that the number of
edges between any two disjoint subsets of $G$ of size $s \leq n/4$ is bounded by
$s\big(np/4+\sqrt{2np \log n}\big)$. Let $H$ be a subgraph of $G$ with maximum degree at
most $np/2-8\sqrt{np \log n}$ and let $G'=G-H$.  We will show that $G'$ has a perfect
matching. This together with the previous lemma immediately implies the assertion (i) of
Theorem \ref{matching}.

By our assumptions on $G$ and $H$ the minimum degree of $G'$ is at least $np/2+6\sqrt{np
\log n}$.
Consider a random partition of the set of vertices
into two equal parts of size $n/2$ each and let $G''$ be the bipartite graph
consisting of all edges of $G'$ intersecting both parts. For every
vertex $v$ its degree in $G''$ is a random variable which has {\em hypergeometric
distribution} with expectation
$d_{G'}(v)/2\geq np/4+3\sqrt{np \log n}=\lambda$. As $p \geq 400 \log n/n$ we have that $\lambda <
5np/12$. Set $t=\sqrt{np \log n}$. Using Chernoff-type estimates
(which are valid also for hypergeometric distributions, see Theorem 2.10 in 
\cite{JanLucRuc00}) we obtain that
$$\mathbb{P}\Big[d_{G''}(v) \leq np/4+2\sqrt{np \log n}=\lambda-t\Big] \leq
e^{-\frac{t^2}{2\lambda}}<n^{-6/5}\ll n^{-1}.$$
Therefore there exist a partition $(V_1, V_2)$ such that the degree of every vertex of the
corresponding bipartite graph
$G''$ is at least $np/4+2\sqrt{np \log n}$. Fix such a partition. We claim
that $G'' \subset G'$  satisfies Hall's condition and therefore has a perfect matching.

We need to check that $|N_{G''}(S)|\geq |S|$ for all $S \subset V_i, i=1,2$.
Let $S' \subset V_i, |S'|>n/4$ with $|N_{G''}(S')|< |S'|$.
If $|N_{G''}(S')| \geq n/4$ let $S=V_{3-i}-N_{G''}(S')$, otherwise let $S$ be any subset of 
$V_{3-i}-N_{G''}(S')$ of size $n/4$.
Note that in both cases $|S| \leq n/4$. Since $N_{G''}(S) \subseteq V_i-S'$, we have 
$$|N_{G''}(S)| \leq |V_i|-|S'|< \min\big(n/4, |V_{3-i}|-|N_{G''}(S')|\big) = |S|.$$
This shows that it is enough to verify inequality $|N_{G''}(S)|\geq |S|$ only for subsets $S$ of
size at most $n/4$. Let $S \subset V_i$ of size $s\leq n/4$.
Note that, by our assumptions there are at most
$s\big(np/4+\sqrt{2np \log n}\big)$ edges of $G$ between any two disjoint sets of size $s$.
On the other hand, it follows from the minimum degree estimate that there are at least
$s\big(np/4+2\sqrt{np \log n}\big)$ edges of $G'' \subset G$ between $S$ and $N_{G''}(S)$.
This implies that $|N_{G''}(S)|\geq s$ and completes the proof of the first part.

\vspace{0.15cm}
\noindent
(ii)\, Partition the vertices of $G=G(n,p)$ into two parts $X$ and $Y$ with sizes
$n/2+1$ and $n/2-1$ respectively. Let
$H$ be a subgraph of $G$ induced by $X$. By part (ii) of Lemma \ref{properties} we have that
a.s. the maximum degree of $H$ is at most $np/2+2\sqrt{np \log n}$. Thus it is enough to
show that $G'=G-H$ has no perfect matching.  Note that the vertices of $X$
form an independent set in $G'$ and therefore should be matched with vertices in $Y$.
On the other hand this is impossible as $|Y|<|X|$.
\hfill $\Box$

\section{Hamilton cycles}
In this section we study the local resilience  of random and
pseudo-random graphs with respect to hamiltonicity. For
illustrative purposes we start by focusing on the case of dense
pseudo-random graphs, since in this case the treatment is
considerably less technical.

\subsection{$\epsilon$-regular graphs}
A graph $G=(V,E)$ on $n$ vertices is called {\em $(d,\epsilon)$-regular} if its minimum
degree is at least $dn$, and for every pair of disjoint subsets $S,T \subset V$ of
cardinalities at least $\epsilon n$, the number of edges between $S$ and $T$ satisfies
$$\left| \frac{e(S,T)}{|S||T|}-d\right| \leq \epsilon.$$

\begin{theo}
\label{pseudo-random}
Let $G$ be a $(d,\epsilon)$-regular graph on $n$
vertices, where $0<\epsilon \ll d\leq 1$ are constants and $n$ is sufficiently large
and let $c_d=8/d$. Then for every subgraph $H$ of $G$ with maximum degree
$\Delta(H) \leq (d/2-c_d\epsilon)n$ the graph $G'=G-H$ contains a Hamilton cycle.
\end{theo}

Using Chernoff-type estimates one can easily check that for every
constant $0<p<1$ and $\epsilon>0$
the random graph $G(n,p)$ is almost surely $(p,\epsilon)$-regular. In particular, we
obtain that by deleting at most $(1/2-\epsilon)np$ edges at every vertex of $G(n,p)$ we
cannot destroy all  Hamilton cycles in this graph. On the other hand, the construction
in part (ii) of Theorem \ref{matching} shows that the constant $1/2$ cannot be improved.

\begin{lemma}
\label{estimate}
Let $G$ be a $(d,\epsilon)$-regular graph on $n$
vertices with $0<2\epsilon<d\leq 1$ and let $H$ be a subgraph of $G$ with maximum degree
$\Delta(H) \leq (d/2-c\epsilon)n$ and $c>2$. Then the graph $G'=G-H$ is
connected and for every $X \subseteq V(G')$ of size at least $\epsilon n$ we have that
$|N_{G'}(X)| \geq (1/2+(c-2)\epsilon)n$.
\end{lemma}

\noindent
{\bf Proof.}\,
Let $X$ be a subset of $V(G')$ of size at least $\epsilon n$ and let $Z$ contains
precisely $\epsilon n$ vertices of $X$.
Consider the set $Y$ of all vertices of $G'-Z$ which have no neighbors in $Z$.
Then $|N_{G'}(X)|\geq |N_{G'}(Z)| \geq n-|Y|-|Z|$. Since $G$ is a $(d,\epsilon)$-regular graph,
by definition, we have $e_G(Z,Y) \geq (d-\epsilon)|Z||Y|$.
On the other hand we know that $e_{G'}(Z,Y) =0$.
This implies that all the edges of $G$ between $Z$ and $Y$ were deleted when we removed $H$.
Since the maximum degree of $H$ is at most $(d/2-c\epsilon)n$, there are at most
$(d/2-c\epsilon)n|Z|$ edges of $H$ incident with vertices in $Z$. Combining these facts
we obtain that
$$(d/2-c\epsilon)n|Z| \geq e_G(Z,Y) \geq (d-\epsilon)|Z||Y|.$$
This implies that
$$|Y| \leq \frac{d/2-c\epsilon}{d-\epsilon} n \leq \left(\frac{1}{2}-(c-1) \epsilon
\right)n,$$
and hence $|N_{G'}(X)|\geq n-|Y|-|Z| \geq (1/2+(c-2) \epsilon)n$.

Next we use this inequality to show that $G'=G-H$ is connected.
Suppose that this is not true and let $X$ be
the vertex set of the smallest connected component of $G'$. Then
$|X| \leq n/2$ and $N_{G'}(X)=X$. Since the minimum degree of $G'$ is at
least $dn-\Delta(H)\geq dn/2$, we have that $|X| \geq dn/2 \geq \epsilon n$
and therefore $|N_{G'}(X)|>n/2\geq |X|$. This contradiction implies that
$G'$ is connected and completes the proof.
\hfill $\Box$

\vspace{0.2cm} \noindent {\bf Proof of Theorem
\ref{pseudo-random}.}\, Let $c_d=8/d$ and let $G'$ be the graph
obtained from $(d,\epsilon)$-regular graph $G$ of order $n$ by
deleting a subgraph $H$ with maximum degree $\Delta(H) \leq
(d/2-c_d\epsilon)n$. To prove the theorem we show that some path
$P$ of a maximal length in $G'$ can be closed to a cycle. As $G'$
is connected by Lemma \ref{estimate}, any non-Hamilton cycle can
be extended to a path covering some additional vertices. Therefore
the assumption about the maximality of $P$ implies that $P$ is a
Hamilton path, and thus the above created cycle is Hamilton as
well. Our approach relies on a variant of the so called
rotation-extension technique, invented by P\'osa in \cite{Pos76}
and applied in several subsequent papers on hamiltonicity of
random and pseudo-random graphs (see, e.g., \cite{KomSze83},
\cite{BolFenFri87}, \cite{KS}).

Let $P=(v_1,v_2,\ldots,v_l)$ be a longest path in $G'$. If $1\le i<l$
and $(v_i,v_l)\in E(G')$, then the path
$P'=(v_1,v_2,\ldots,v_i,v_l,v_{l-1},\ldots,v_{i+1})$ is also of maximal
length. We say that $P'$ is a {\em rotation} of $P$ with {\em fixed
endpoint} $v_1$, {\em pivot} $v_i$ and {\em broken edge}
$(v_i,v_{i+1})$ (the reason for the last term being the fact that
$(v_i,v_{i+1})$ is deleted from the edge set of $P$ to get $P'$). We can
then rotate $P'$ in a similar fashion to get a new path $P''$ of the
same length, and so on. A subset $I$ of consecutive points of $P$
is called a {\em segment}. Given set a $S \subseteq P$, a vertex $v \in S$ is an {\em
interior point} of $S$ with respect to $P$ if both neighbors of $v$ along $P$ lie in $S$.
The set of all interior points of $S$ will be denoted by $int(S)$.

Let $t=d/(5\epsilon)$ and let $I_1, I_2, \ldots, I_t$ be a
partition of path $P$ into $t$ segments of length $|P|/t \leq n/t=
(5/d)\epsilon n$. By definition, all the neighbors of $v_1$ and
$v_l$ in $G'$ belong to $P$. Since the minimum degree of $G'$ is
at least $dn-\Delta(H)\geq dn/2$, there exist a segment $I_p$
which contains at least $(dn/2)/t \geq 5\epsilon n/2$ neighbors of
$v_1$. Similarly there exist a segment $I_q$ which contains at
least $5\epsilon n/2$ neighbors of $v_l$. If $p=q$ divide $I_p$
into two segments $J_1, J_2$ such that both contains exactly
$5\epsilon n/4$ neighbors of $v_1$. Clearly, one of the segments
$J_i$ contains at least $5\epsilon n/4$ neighbors of $v_l$.
Without loss of generality suppose that it is $J_2$. Hence we
obtain two disjoint segments of $P$ such that $|J_1|, |J_2| \leq
(5/d)\epsilon n$, $int(J_1)$ contains at least $5\epsilon
n/4-2>\epsilon n$ neighbors of $v_1$ and $int(J_2)$ contains at
least $\epsilon n$ neighbors of $v_l$. In case $p\not =q$ simply
take $J_1=I_p$ and $J_2=I_q$.

For every neighbor $u$ of $v_1$ in $int(J_1)$ rotate $P$ using $u$
as pivot and keeping $v_l$ fixed. Let $A$ be the set of endpoints
of paths obtained by such rotations. For every $a\in A$ let $P(a,
v_l)$ be the corresponding new maximum length path with endpoints
$a$ and $v_l$ directed from $a$ to $v_l$. By the above discussion,
$A$ is subset of $J_1$ of size at least $\epsilon n$ and all the
edges of $P$ outside $J_1$ are unbroken and belong to every path
$P(a, v_l)$. Moreover, each $P(a, v_l)$ traverse the segments of
$P$ outside $J_1$ in exactly  the same direction.

Now for every $a \in A$ and for every neighbor $w$ of $v_l$ in
$int(J_2)$ rotate $P(a, v_l)$ using $w$ as pivot and keeping $a$
fixed. For every $a \in A$, let $B(a)$ be the set of endpoints of
paths obtained by such rotations. Note that, since $J_2$ is
disjoint from $J_1$, all the paths $P(a, v_l)$ traverse the
segment $J_2$ in the same direction. This implies that the set
$B(a)$ does not depend on $a$. Call this set $B$. The size of $B$ is
at least $\epsilon n$. Furthermore, $B$ is completely contained in
$J_2$. For every $a \in A$ and $b \in B$, let $P(a, b)$ be the
corresponding new maximum length path with endpoints $a$ and $b$,
directed from $a$ to $b$. Again it is easy to see that in all the
paths $P(a, b)$ the edges outside $J_1 \cup J_2$ are unbroken and
these paths traverse the segments of $P$ outside $J_1 \cup J_2$ in
exactly the same direction.

Consider a point $v$  in the interior of $P-(J_1 \cup J_2)$
adjacent to some $b \in B$. For every $a \in A$ and for every such
$v$ we can rotate $P(a, b)$ using $v$ as pivot and keeping $a$
fixed. Since all the paths $P(a, b)$ traverse the edges of $P-(J_1
\cup J_2)$ in the same direction,  we conclude that the set of
endpoints of paths obtained by such rotations does not depend on
$a$. Denote this set by $C$. Then for every $a \in A$ and $c \in
C$ there exists a maximum length path with endpoints $a$ and $c$.
Note that the size of $C$ is at least
$|N_{G'}(B)|-|J_1|-|J_2|-O(1)$. Furthermore $|B| \geq  \epsilon
n$. Therefore, by Lemma \ref{estimate} we have that
\begin{eqnarray*}
|C| &\geq&  \left(\frac{1}{2}+(c_d-2) \epsilon\right)n-2\big(5/d\big)\epsilon n-
O(1)\\
&>&\left(\frac{1}{2}+\big(8/d-2\big) \epsilon\right)n -\big(11/d\big)\epsilon n\\
&\geq& \left(\frac{1}{2} -\big(5/d\big) \epsilon\right)n.
\end{eqnarray*}
On the other hand, since also $|A| \geq  \epsilon n$ the same lemma implies that
$$|N_{G'}(A)| \geq \left(\frac{1}{2}+(c_d-2) \epsilon\right)n
\geq \left(\frac{1}{2}+\big(6/d\big) \epsilon\right)n
>n-|C|.$$
Thus we conclude that $C \cap N_{G'}(A)$ is nonempty, i.e., there
is an edge connecting $A$ and $C$ and hence closing a maximal path
to a cycle. This completes the proof. \hfill $\Box$

\subsection{Sparse random graphs}

Now we are going to consider the general case when $p$ can be as
small as $\log^4 n/n$.

\begin{theo}
\label{hamilton}
For every fixed $\epsilon>0$ and $p \geq\log^4 n/n$ the random graph $G(n,p)$
almost surely has the following property.
If  $H$ is a subgraph of $G=G(n,p)$ with maximum degree
$\Delta(H) \leq (1/2-\epsilon)np$ then $G'=G-H$ contains a Hamilton cycle.
\end{theo}

The proof of this theorem is much more technical than the previous
one. Let us start with two (rather routine) lemmas about the edge
distribution of $G(n,p)$.

\begin{lemma}
\label{prop-random}
If $p \geq \log^4 n/n$ then almost surely

\noindent $(i)$\, For every subset $A, |A|=a$ of $G(n,p)$ the
number of edges between $A$ and its complement $V(G)-A$ is $(1+o(1))a(n-a)p$.

\noindent $(ii)$\, The number of edges between any two disjoint subsets $A, B$ of
$G(n,p)$ of size $|A|=a$ and $|B|=b \leq \min\big(\frac{anp}{15}, \frac{n}{15}\big)$
is less than $anp/2$.

\noindent $(iii)$\, The number of edges between any two disjoint subsets $A, B$ of $G(n,p)$
of size $|A|=a \geq n/\log^3 n$ and $|B|=b \geq n/\log^{1/4} n$ is at least
$(1+o(1))abp$.
\end{lemma}

\noindent
{\bf Proof.}\,  (i) By symmetry it is enough to prove the claim for $a \leq n/2$.
In this case, using Chernoff-type estimates, we obtain
 $$\mathbb{P}\left[ \Big|e(A, V(G)-A) -a(n-a)p\Big| \geq t= \frac{a(n-a)p}{\log n}\right] \leq
 e^{-\Omega\left(\frac{t^2}{2a(n-a)p}\right)}
 \leq e^{-\Omega\left(\frac{a(n-a)p}{\log^2 n}\right)} \ll e^{-2a\log n}.$$
Therefore the probability that there is set $A, |A| \leq n/2$ which violates condition (i) is
at most
$$\sum_{a \leq n/2} {n \choose a}e^{-2a\log n} \leq \sum_{a \leq n/2} n^a e^{-2a\log n}=
\sum_{a \leq n/2} e^{-a\log n}= o(1).$$

\noindent
(ii)\, The probability that there is a pair of sets $A, B$ that violates this claim is at most
\begin{eqnarray*}
\mathbb{P} &\leq& \sum_{a, b}{n \choose a} {n \choose b} {ab \choose anp/2}p^{anp/2} \leq
 \sum_{a,b}{n \choose a} \left(\frac{en}{b}\right)^b \left(\frac{eab}{anp/2}\right)^{anp/2}
 p^{anp/2}\\
 &\leq& \sum_{a,b}{n \choose a} \left(\frac{en}{b}\right)^{anp/15}
 \left(\frac{2eb}{n}\right)^{anp/2} \leq \sum_{a,b}n^a
\left(2^{1/2}e^{3/5}\big(b/n\big)^{2/5}\right)^{anp}\\
&\leq& \sum_{a,b}n^a \left(\frac{2^3e^3}{15^2}\right)^{anp/5}
 \leq n \sum_a e^{a \log n}\, 0.8^{a\log^4 n/5}= o(1).
\end{eqnarray*}

\noindent
(iii)\, The number of choices for sets $A, B$ is at most $2^{2n}$. On the other hand,
using again Chernoff-type estimates, we have that
$$\mathbb{P}\left[ e(A, B) -abp \leq -t= -\frac{abp}{\log^{1/4} n}\right] \leq
 e^{-\frac{t^2}{2abp}}
 \leq e^{-\frac{abp}{2\log^{1/2} n}} \leq e^{-n \log^{1/4} n/2} \ll 2^{-2n}.$$
This completes the proof of the lemma. \hfill $\Box$

\begin{lemma}
\label{expansion}
For $p \geq \log^4 n/n$ and fixed $\epsilon>0$
let $H$ be a subgraph of $G=G(n,p)$ with maximum degree
$\Delta(H) \leq (1/2-\epsilon)np$. Then
a.s. the graph $G'=G-H$ is connected,  every subset $A$ of $G'$
of size $a \leq n/\log^4 n$ satisfies that $|N_{G'}(A)| \geq a \log^4 n/15$
and every subset $B$ of size $b \geq n/\log^3 n$ has
$|N_{G'}(B)| \geq \frac{1+\epsilon}{2}n$.
\end{lemma}

\noindent
{\bf Proof.}\, Let $H$ be a subgraph of $G=G(n,p)$ with maximum degree
$\Delta(H) \leq (1/2-\epsilon)np$. Consider a  subset $A$ of $G'=G-H$
of size $a \leq n/\log^4 n$. By part (i) of the previous lemma
we know that there are $(1+o(1))a(n-a)p=(1+o(1))anp$ edges of $G(n,p)$
between $A$ and its complement $V(G)-A$. On the other hand there are at most
$(1/2-\epsilon)anp$ edges of $H$ incident with vertices in $A$.
Combining these facts we obtain that
$$e_{G'}(A, V-A) \geq (1+o(1))anp - (1/2-\epsilon)anp \geq anp/2.$$
Let $Y$ be the set of all vertices from $V-A$ which are adjacent in $G'$ to
some vertex in $A$. Then $e_{G}(A,Y)\geq e_{G'}(A,Y) \geq anp/2$.
Hence, by part (ii) of Lemma \ref{prop-random}
we have that
$$|N_{G'}(A)| \geq |Y| \geq \min\left(\frac{anp}{15}, \frac{n}{15}\right) \geq
\frac{\log^4 n}{15} a.$$

Let $B$ be a subset of $G'$ of  size at least $n/\log^3 n$ and let
$D \subseteq B$ contains precisely $|D|=d=n/\log^3 n$ vertices from $B$.
Consider the set $Z$ of all vertices of $G'-D$ which have no neighbors in $D$.
Then $|N_{G'}(B)|\geq N_{G'}(D) \geq n-|D|-|Z|$.
Note that, by definition of $Z$, all the edges of $G(n,p)$ between $D$ and $Z$ were deleted
when we removed $H$.
Since the maximum degree of $H$ is at most $(1/2-\epsilon)np$, there are at most
$(1/2-\epsilon)ndp$ edges of $H$ incident with vertices in $D$. Using this fact
together with part (iii) of Lemma \ref{prop-random}, we obtain
$$ (1+o(1))d|Z|p \leq e_{G}(D,Z) \leq (1/2-\epsilon)ndp,$$
and therefore $|Z| \leq \big(1/2-\epsilon-o(1)\big)n$.
This implies that
$$|N_{G'}(B)|\geq n-|Z|-|D| \geq \big(1/2+\epsilon-o(1)\big)n-
n/\log^3 n >\frac{1+\epsilon}{2}n.$$

By the above discussion we have that $|N_{G'}(C)|>|C|$ for every subset $C, |C| \leq n/2$.
Thus every connected component of $G'$ has size larger than $n/2$
and therefore $G'$ is connected.  \hfill $\Box$

\vspace{0.2cm} With these lemmas in hand, Theorem \ref{hamilton}
is a corollary of  the following deterministic result, which
ensures the existence of a hamiltonian cycle in a
graph provided some information about edge expansion. 

\begin{theorem}
\label{hamilton-cycle}
For every fixed $\epsilon>0$ and sufficiently large $n$, the following holds.
Let $G$ be a connected graph of order $n$
such that every subset $U$ of $G$ of size
at most $n/\log^4 n$ satisfies $|N_{G}(U)| \geq |U| \log^4 n/15$
and every subset $W$ of size at least $n/\log^3 n$ has
$|N_{G}(W)| \geq \frac{1+\epsilon}{2}n$. Then $G$ contains a Hamilton cycle.
\end{theorem}

\noindent
The proof of Theorem \ref{hamilton-cycle} uses roughly the same
approach as described in the previous subsection but is much more involved. It is enough to show that some 
longest path in $G$ can be closed to a cycle. Since $G$ is connected this cycle must be 
hamiltonian. We start with some fixed longest path $P$ and use the 
rotation-extension technique to construct two sets of vertices $A$ and $B$ of nearly linear size such 
that for every vertex $a \in A$ and every vertex $b \in B$
there is a path $P(a, b)$ from $a$ to $b$ of maximum length. Also,
we have a collection of segments of the original path $P$
with total length $|P|-o(n)$ such that all paths  $P(a,b)$ contain these segments untouched.
Moreover, each of these segments is traversed by all paths $P(a,b)$ in the same direction.
Once we have such a collection of maximum paths $P(a,b)$, we perform one additional round
of rotations to find a maximum path which can be closed to a cycle. This part of the argument is 
exactly the same as in the proof of Theorem \ref{pseudo-random}.

To construct sets $A, B$ as above in the case of dense graphs, it was enough
to rotate each endpoint of the original path 
$P$ once. Unfortunately, in the case of sparse graphs, after one round of rotations one can only obtain sets 
$A, B$ of polylogarithmic size. The goal of the proof is to amplify the size of $A$ and $B$ to be close to linear.
This is achieved in roughly $\frac{\log n}{\log \log n}$ rounds of rotations, where in each round we 
first rotate the left endpoint of all maximum paths $P(a,b)$ and then the right endpoint.
At the end of each round, we obtain new sets $A$ and $B$, which are larger than the previous ones
by a factor of $\log n $.

\vspace{0.25cm} \noindent {\bf Proof.}\, 
We use the same notation as in the proof of 
Theorem \ref{pseudo-random}. Let $P=(v_1, v_2, \ldots,
v_k)$ be the longest path in $G$. We fix a canonical direction on
the path $P$ from $v_1$ to $v_k$, and call $v_1$ the left end and
$v_k$ the right end.  If $1\le i<k$
and $(v_i,v_k)\in E(G)$, then the path
$P'=(v_1,v_2,\ldots,v_i,v_k,v_{k-1},\ldots,v_{i+1})$ is also of maximal
length. We say that $P'$ is a {\em rotation} of $P$ with {\em fixed
endpoint} $v_k$, {\em pivot} $v_i$ and {\em broken edge}
$(v_i,v_{i+1})$ (the reason for the last term being the fact that
$(v_i,v_{i+1})$ is deleted from the edge set of $P$ to get $P'$). We can
then rotate $P'$ in a similar fashion to get a new path $P''$ of the
same length, and so on. We only perform rotations whose broken edges are the edges of the original path $P$.
For example, if we perform a rotation with pivot $v_i$ we will not use vertices $v_{i-1}, v_{i+1}$ as pivots for future rotations.
A subset $I$ of consecutive points of $P$
is called a {\em segment}. Given a set $S \subseteq P$, a vertex $v \in S$ is an {\em
interior point} of $S$ with respect to $P$ if both neighbors of $v$ along $P$ lie in $S$.
The set of all interior points of $S$ will be denoted by $int(S)$.

The proof consists of two parts. First, for every $1 \leq t \leq r= \frac{\log n}{\log \log n}-5$, we use 
induction to construct
segments $I_1, \ldots, I_t$ and $J_1, \ldots, J_t$
of the original path $P$ and two subsets $A_t, B_t$ of size $\log^{t+2} n$ with the following
properties.

\vspace{0.15cm}
\noindent
{\bf Properties:}
\begin{enumerate}
\item All segments $I_{\ell}, J_{\ell'}$ have length at most
$O(n/\log n)$ for every $1 \leq \ell, \ell' \leq t$. Segments
$I_{\ell}$ are disjoint from all segments $J_{\ell}$. For every two
indices $\ell \not = \ell'$, the corresponding segments $I_{\ell}$ 
and $I_{\ell'}$ are either disjoint or the same. Similar conditions
hold for segments  $J_{\ell}, J_{\ell'}$. Also, $A_t
\subset I_t$ and $B_t \subset J_t$.

\item For every $a \in A_t$ and $b \in B_t$, we have a path of
maximum length going  from $a$ to $b$. This path was obtained from
$P$ by $t$ iterations of rotations, where in each iteration we
first rotate  the left end point and then the right end point of a
path obtained in the previous iteration. Therefore, this path has $2t$ broken
edges. The edge which gets broken in the $i$-th iteration when we 
rotate the left end point is denoted by $e_i(a)$. Similarly,
$e_i(b)$ is the edge which gets broken in the $i$-th iteration
when we rotate  the right end point. These $2t$ edges partition
$P$ into $2t+1$ segments $Q_1(a,b), \ldots, Q_{2t+1}(a,b)$, where
the indices $1,2 \dots, 2t+1$ correspond to the order in which the
segments appear on the path $P(a,b)$ when we go from $a$ to  $b$. 

\item For any $a \in A_t$, the corresponding edges $e_1(a), \ldots,
e_t(a)$ depend only on $a$, regardless of the right end point $ b \in
B_t$. Also, the edge $e_{\ell}(a)$ belongs to $I_{\ell}$ for all $a
\in A_t$ and $1\leq \ell \leq t$. Similarly, the edges $e_1(b),
\ldots, e_t(b)$ depend only on $b$, and $e_{\ell}(b)$ belongs to $
J_{\ell}$ for all $b \in B_t$ and $1\leq \ell \leq t$. Moreover, 
the order in which the $2t$ edges $\{e_{\ell}(a)\}
\cup\{e_{\ell}(b)\}$ appear on the path $P$ is universal (i.e.,
it does not depend on the pair $(a,b)$).

\item There is a universal pair $(\pi, z)$
(i.e., it does not depend on the endpoints $a,b$), where
$\pi$ is a permutation in $S_{2t+1}$ and $z$ is a binary vector of
length $2t+1$. Permutation $\pi$ records  the order in which the  
segments $Q_i(a,b)$ appear on the path $P$. For example $\pi(i)=j$
means that $Q_i(a,b)$ is the $j$-th segment from the left on $P$. 
The vector $z$ records the direction in which these segments are  
traversed by $P(a,b)$: $z_i=0$ if  $P(a,b)$ traverses $Q_i(a,b)$  
from left to right (in the direction of $P$) and $z_i=1$
otherwise. Therefore, although segments $Q_i(a,b)$ do depend on the endpoints $a,b$, we
have that these segments appear on paths $P$ in exactly the same order, and
also the direction in which each segment $Q_i(a,b)$ is traversed 
is the same for all paths $P(a,b)$.
\end{enumerate}

Suppose that we have constructed sets $A_r, B_r$ of size 
$\log^{r+2} n$ and segments $I_{1}, \ldots, I_{r}$ and 
$J_{1}, \ldots, J_{r}$ which satisfy properties 1--4 for $r=\frac{\log n}{\log \log n}-5$.
Then we can complete the proof of the theorem as follows.     
Note that both $A_r$ and $B_r$ have size $\log^{r+2} n=n/\log^3
n$, and all paths $P(a,b)$ have the same vertex set as path $P$. Also by maximality 
of path $P(a,b)$, all neighbors of $a$ also belong to this path. Thus 
the neighbors of all vertices in $A_r$ belong to path $P$.
Let $I=\cup_{\ell=1}^rI_{\ell}$, $J=\cup_{\ell=1}^rJ_{\ell}$, and let $v$ be a
point in the interior of $P-(I \cup J)$  adjacent to some $a \in
A_r$. By the above discussion, the number of such vertices $v$ is at least
$|N_{G}(A_r)|-|I|-|J|-O(r)$.
For every $b \in B_r$ and for every such $v$, we can rotate
$P(a, b)$ using $v$ as pivot and keeping $b$ fixed. From
properties 1--4, it is easy to check that for every edge of $P-(I \cup J)$, 
there is an index $1 \leq i \leq 2r+1$ such that for all paths $P(a,b)$ this edge belongs 
to the segment $Q_i(a,b)$. Therefore, every edge of $P-(I \cup J)$ is traversed 
in the same direction by all paths $P(a, b)$.
Thus, the set of endpoints of new
paths obtained by such rotations does not depend on $b$. Denote this
set by $C$. For every $b \in B_r$ and $c \in C$ there exists a
maximum length path with endpoints $c$ and $b$. As we already mentioned, the
size of $C$ is at least $|N_{G}(A_r)|-|I|-|J|-O(r)$ and
$|A_r|=n/\log ^3 n$. By our  assumptions, we have that
$$|C| \geq  \frac{1+\epsilon}{2}n-2r\, \frac{n}{\log n}-O(r) >\frac{1+\epsilon}{2}n
-o(n)>n/2.$$ On the other hand, since $|B_r|= n/\log ^3 n$, the   
expansion assumption yields
$$|N_{G}(B_r)| \geq \frac{1+\epsilon}{2}n  >n-|C|.$$
This implies that  $C \cap N_{G}(B_r)$ is not empty. So, there is
an edge connecting $B_r$ and $C$ which closes a maximal path to a
cycle. The proof is complete.  \hfill $\Box$

\vspace{0.15cm}
Next we show, using induction, how to construct segments 
$I_{1}, \ldots, I_{t}$ and $J_{1}, \ldots, J_{t}$, 
and sets $A_t, B_t$ of size $\log^{t+2} n$, satisfying properties 1--4.
We start with a detailed description of the basis case when $t=1$.

\vspace{0.15cm}
\noindent
{\bf Construction for $t=1$.}\,  By our assumption, the minimum degree $\delta(G)$ (of $G$) is at least 
$\log^4 n/15$. Due to the maximality of path $P$, all the neighbors of $v_1$ and 
$v_k$ belong to $P$. Set $s=\log n/40$ and partition $P$ into $s$ segments 
of length $|P|/s \leq n/s=O(n/\log n)$. By the pigeonhole principle, 
there are two segments $S_1$ and $S_2$ containing at least $\delta(G)/s 
\geq 8\log^3 n/3$ neighbors of $v_1$ and $v_k$, respectively. If $S_1\not 
=S_2$, let $I_{1}=S_1$ and $J_{1}=S_2$. Otherwise, divide $S_1$ into two 
segments, each containing at least $4\log^3 n/3$ neighbors of $v_1$. One of 
these segments must contain at least $4\log^3 n/3$ neighbors of $v_k$. 
Call this segment $J_{1}$ and the other one $I_{1}$. In both cases, we obtain 
two disjoint segments $I_{1}, J_{1}$ of length $O(n/\log n)$ such that 
$int(I_{1})$ contains at least $\log^3 n$ neighbors of $v_1$ and $int(J_{1})$ 
contains at least $\log^3 n$ neighbors of $v_k$.

For every neighbor $u$ of $v_1$ in $int(I_{1})$, rotate $P$ using $u$
as pivot and keeping $v_k$ fixed. Let $A_1$ be the set of
endpoints of paths obtained by such rotations.  For $a\in A_1$,
let $P(a, v_k)$ be the corresponding new maximum  path with
endpoints $a$ and $v_k$, directed from $a$ to $v_k$. Also, for every
$a \in A_1$ and for every neighbor $w$ of $v_k$ in $int(J_{1})$,
rotate $P(a, v_k)$ using $w$ as pivot and keeping $a$ fixed. Let
$B_1$ be the set of endpoints of paths obtained by such rotations.
Since $J_{1}$ is disjoint from $I_{1}$, all paths $P(a, v_k)$ traverse
the segment $J_{1}$ in the same direction. This implies that the set
$B_1$ does not depend on $a$.  Both $A_1$ and $B_1$ are of size at least
$\log^3 n$. By deleting some vertices, we can assume that
$|A_1|=|B_1|=\log^3 n$ for convenience. We also know that $A_1
\subset I_{1}$, $B_1 \subset J_{1}$, and for every $a \in A_1$ and $b
\in B_1$ we have a maximum length path $P(a,b)$ with endpoints $a$
and $b$. 

Each path $P(a,b)$ was obtained by two rotations, one
from the left and one from the right. Furthermore, it has two
broken edges $e_1(a)$ (obtained when we rotated $v_1$) and
$e_1(b)$ (obtained when we rotated $v_k$). Note that the edge
$e_1(a)$ belongs to $I_{1}$ and does not depend on $b$. Similarly, the edge $e_1(b)$ belongs to
$J_{1}$ and does not depend on $a$. In particular, the order in which $e_1(a)$ and $e_1(b)$
appear on the path $P$ is the same order in which the segments $I_{1}$ and $J_{1}$ appear.
Therefore this order is universal  and does not depend on the
endpoints $a$ and $b$. The edges $e_1(a)$ and $e_1(b)$ partition the  path $P$ into 3
segments $Q_1(a,b), Q_2(a,b), Q_3(a,b)$, where the indices $1,2,3$
correspond to the order these segments appear on  $P(a,b)$ when we
go from $a$ to $b$. We associate with each path $P(a,b)$ a 
permutation in $S_3$ (the permutation
group on $1,2,3$) and a binary vector of length $3$. The
permutation records the order in which the segments
$Q_i(a,b)$ appear on the path $P$. The vector records the
direction in which these segments were traversed by the path
$P(a,b)$ as follows. The coordinate $i$ of the vector is $0$ if  $P(a,b)$ traverses
$Q_i(a,b)$ from left to right (in the direction of $P$) and
is $1$ otherwise. From our construction it follows that, although
the location of the segments $Q_i(a,b)$ does depend on the actual
endpoints $a$ and $ b$, there is a universal pair $(\pi, z)$ such that all segments
$Q_i(a,b)$ appear on path $P$ in the order of $\pi$, and all paths $P(a,b)$ traverse the 
corresponding segments $Q_i(a,b)$ in the direction defined by the $z_i$. In
particular, every edge outside $I_{1,1} \cup J_{1,1}$ is unbroken, has the
same index $1\leq i \leq 3$ of interval $Q_i(a,b)$ to which it
belongs for all paths $P(a,b)$, and therefore is traversed in the
same direction by all these paths. \hfill $\Box$

\vspace{0.15cm} \noindent {\bf Induction step:}\, Suppose we have segments 
$I_1, \ldots, I_t$ and $J_1, \ldots, J_t$ of the path $P$ and two subsets 
$A_t, B_t$ of size $\log^{t+2} n$, satisfying properties 1--4. We construct new sets of end points 
$A_{t+1}$ and $B_{t+1}$ by first rotating the left end points of all paths 
$P(a,b), a \in A_t, b \in B_t$, keeping right end points fixed, and then 
vice versa. When we rotate the left end points of all paths $P(a,b)$, we 
obtain new segments $I_1, \ldots, I_t, I_{t+1}$ and $J_1, \ldots, J_t$, and 
sets $A_{t+1}$ and $B_t$. The set $A_{t+1}$ will be larger than the original set $A_t$
by a factor $\log^2 n$.
On the other hand, the new set $B_t$ 
might decrease by a factor $2(t+1)$ during this operation. Note that we 
slightly abusing notation here, since both segments $I_1, \ldots, I_t$, 
$J_1, \ldots, J_t$ and also set $B_t$ are changing during the rotation, and 
are not what they were originally. Next we rotate right endpoints of all 
paths $P(a,b), a \in A_{t+1}, b \in B_t$ (here $B_t$ is the new set we got 
after rotating left end points). In the end of the rotation from the right 
we obtain new segments $I_1, \ldots, I_t, I_{t+1}$ and $J_1, \ldots, 
J_t, J_{t+1}$, and sets $A_{t+1}$ and $B_{t+1}$. After the second round of 
rotations the size of the set $B_{t+1}$ will be larger than the set $B_t$
by a factor $\log^2 n$, but the set $A_{t+1}$ might shrink by a factor 
of $2(t+1)$. Overall, after rotating left and right end points,
the sizes of the new sets $A_{t+1}$, $B_{t+1}$ are larger then the sizes of the original sets
$A_t, B_t$ by a factor of at least $\Omega(\log^2 n/t)\gg \log n$. 
Moreover 
we perform rotations so that we preserve all the properties 1--4. Thus, at 
the end of the induction step we obtain new sets $A_{t+1}, B_{t+1}$ of 
size $\log^{t+3} n$ (this is in fact a lower bound, but we can always 
delete extra vertices) and two collections of segments $I_1, \ldots, 
I_{t+1}$ and $J_1, \ldots, J_{t+1}$ satisfying properties 1--4.

Now we give a detailed description of how we rotate the left end points of 
paths $P(a,b)$. Rotation of the right end points is done similarly and we 
will omit its description here. By the expansion properties of $G$ we have 
that $|N(A_t)| \geq (\log^4 n)|A_t|/15$. Due to the maximality of paths 
$P(a,b)$ and the fact that they all have the same vertex set as $P$, 
all neighbors of vertices in $A_t$ belong to $P$. To give
further details of the rotation procedure we need to consider several 
cases. The first case is when there are at least $|N(A_t)|/3$  neighbors of $A_t$  outside the 
set $\cup_{\ell}(I_{\ell}\cup J_{\ell})$. This a basic case which is also used to analyze the other two cases.
The second case is when there are at least 
$|N(A_t)|/3$ neighbors of $A_t$ inside $\cup_{\ell}J_{\ell}$. The third (final) case is when there
are at least $|N(A_t)|/3$ neighbors of $A_t$ inside $\cup_{\ell}I_{\ell}$. 
Clearly one of these three cases always happens.

{\bf Case 1.}\, Suppose that there are at least $(\log^4 n)|A_t|/45$ vertices of
$N(A_t)$ which do not belong to any segment $I_{\ell}$ or
$J_{\ell}$. Since $P\setminus \cup_{\ell} (I_{\ell} \cup
J_{\ell})$ is a union of at most $2t+1$ segments, we can take one
of them, which we call $S$, that contains at least $\Omega\big((\log^4
n)|A_t|/t\big)$ vertices from $N(A_t)$. Partition $S$ into $\log n$
segments of length $|S|/\log n \leq n/\log n$. Then the
interior of one of these segments contains at least $\Omega\big((\log^3
n)|A_t|/t\big) \gg (\log^2 n) |A_t|$ points from $N(A_t)$. Denote this
segment by $I_{t+1}$.

Let $P(a,b)$ be an arbitrary maximum  path with $a \in A_t, b \in
B_t$. This path has $2t$ broken edges $\{e_{\ell}(a)\}
\cup\{e_{\ell}(b)\}$ which partition it into $2t+1$ segments
$Q_1(a,b), \ldots, Q_{2t+1}(a,b)$, which are also segments of the
original path $P$. As usual, the  indices of these segments
correspond to the order in which they appear on $P(a,b)$ when we
go from $a$ to $b$. Since $e_{\ell}(a) \in I_{\ell}, e_{\ell}(b)
\in J_{\ell}$ and $I_{t+1}$ is located outside $\cup_{\ell}
(I_{\ell} \cup J_{\ell})$, we have that $I_{t+1}$ contains none of
these edges and its position with respect to $\{e_{\ell}(a)\}
\cup\{e_{\ell}(b)\}$ is exactly the same for all $a \in A_t, b \in
B_t$. Therefore there is an integer $ 1 \leq q \leq 2t+1$, which
does not depend on the endpoints $a,b$, such that $I_{t+1}$ is a
subset of the $q$-th segment (from the left) of $P\setminus
\{e_{\ell}(a)\} \cup\{e_{\ell}(b)\}$. Let $1 \leq h \leq 2t+1$ be
such that $\pi(h)=q$. By the definition of the permutation $\pi$,
 $I_{t+1}$ is a subset of $Q_{h}(a,b)$ for all paths $P(a,b)$.
Without loss of generality, suppose that $z_h=0$, i.e., all paths
$P(a,b)$ traverse the segment $Q_{h}(a,b)$ from left to right (in
the direction of path $P$). The case when $z_h=1$ can be treated
similarly.

Let $v$ be a neighbor of $a$ in the interior of $I_{t+1}$ and let
$(a',v)$ be the edge of $P$ (and hence edge of $P(a,b)$) where
$a'$ is the vertex which is immediately to the left of $v$. Then we
can rotate the path $P(a,b)$, using $v$ as a pivot and keeping $b$
fixed, to obtain  a new path $P(a',b)$. Put $a' \in A_{t+1}$, and
note that for all $b \in B_t$, we will have the same new broken
edge $(a',v)$, which we denote by $e_{t+1}(a')$. Also, for all $1
\leq \ell \leq t$, define $e_{\ell}(a')=e_{\ell}(a)$. By our
construction, it is easy to see that for all $a'$ and $b$, the
order in which the $2t+1$ edges $\{e_{\ell}(a')\} \cup\{e_{\ell}(b)\}$
appear on the path $P$ is universal (i.e, does not depend on the
pair $(a',b)$).

The edges $\{e_{\ell}(a')\} \cup\{e_{\ell}(b)\}$ partition $P$
into $2t+2$ segments $Q'_1(a,b), \ldots, Q'_{2t+2}(a,b)$, where
the indices correspond to the order in which these segments are
traversed by the path $P(a',b)$. New segments are obtained from
segments $Q_1(a,b), \ldots, Q_{2t+1}(a,b)$ as follows. $Q'_1(a,b)$
is the left part of  $Q_{h}(a,b)$ from the beginning of $Q_{h}(a,b)$
to $a'$, and $Q'_{h+1}(a,b)$ is the right part of  $Q_{h}(a,b)$ from
the pivot vertex $v$ to the end. For all $ 2\leq p \leq h$, we
have that $Q'_p(a,b)=Q_{h-p+1}(a,b)$. Moreover, for all $ h< p
\leq 2t+2$,  $Q'_p(a,b)=Q_{p+1}(a,b)$. Thanks to the
correspondence between the old and new segments,  the order in
which the segments $\{Q'_p(a,b)\}$ appear on the path $P$ does not
depend on the path $P(a',b)$. Therefore, the permutation $\pi'$
which records this order is the same for all new paths. Finally, we
will show that the new vector $z'$ that records the directions in
which $P(a',b)$ traverses the segments $\{Q'_p(a,b)\}$ is also
universal. Note that after the rotation of $P(a,b)$ with fixed $b$
and pivot $v$, the new path has opposite direction on the part of
$P(a,b)$ before $v$, and same direction as $P(a,b)$ on the part after $v$.
Therefore, we have that $z_1=1-z_h$ and $z_p=1-z_{h-p+1}$ for all $2
\leq p \leq h$, and $z_p=z_{p+1}$ for all $h \leq p \leq 2t+2$. Now
we can conclude that the new family of maximum paths possesses 
properties 1--4.

\vspace{0.15cm} {\bf Case 2.}\, Suppose that there are at least $(\log^4 n)|A_t|/45$ vertices of
$N(A_t)$ inside $\cup_{\ell} J_{\ell}$. Then there is
an index $\ell^*$ such that $J_{\ell^*}$ contains at least
$\Omega\big((\log^4 n)|A_t|/t\big)>(\log^3 n)|A_t|$ vertices from
$N(A_t)$. Since some of the segments $J_{\ell}$ might be equal,
let $L$ be the set of all indices $1 \leq \ell \leq t$ such that
$J_{\ell}=J_{\ell^*}$. Partition $J_{\ell^*}$ into $2t$ segments
$S_1, \ldots, S_{2t}$, each containing $(\log^3 n)|A_t|/(2t)\gg
(\log^2 n)|A_t|$ vertices of $N(A_t)$. Since for every fixed $b \in
B_t$, the number of broken edges $\{e_\ell(b)\}$ is at most $t$, at
least half of all the segments $S_j$ contain no such edge. By
averaging,  there is a segment $S_i$ such that for at least half
of the vertices $b \in B_t$ this segment contains no edges
$\{e_\ell(b)\}$. Let $B'_t$ be the collection of those $b$ where
$\{e_\ell(b)\}$ is not in $S_i$.  For every $b \in B'_t$, consider the
location of the segment $S_i$ with respect to the edges
$\{e_\ell(b), \ell \in L\}$. Note that by properties 1--4, these
edges appear in the same order inside $J_{\ell^*}$ and partition
it into at most $t+1$ segments, one of which contains $S_i$.
Therefore, there is a subset $B''_t \subset B'_t$ of size at least
$|B'_t|/(t+1)$ such that for all $b \in B''_t$, the segment $S_i$
has exactly the same location with respect to the  edges
$\{e_\ell(b)\}$, $\ell \in L$.  $J_{\ell^*}\setminus S_i$ consists
of two disjoint segments: $J'$ on the left  and $J''$ on the
right. Moreover, let $L'$ be the set of indices $\ell \in L$ such
that $e_\ell(b)$ appears on the left of $S_i$, and let $L''$ be the set
of indices $\ell \in L$ such that $e_\ell(b)$ appear on the right
of $S_i$.

Define new segments $J_\ell=J'$ for $\ell \in L'$ and $J_\ell=J''$
for $\ell \in L''$. Redefine $B_t=B''_t$. Note that since we only
delete points from $B_t$, the  new segments and the sets $A_t$ and
$B_t$ still satisfy all properties 1--4. On the other hand, we now
have the segment $S_i$ of length at most $n/ \log n$ (since it is
a subset of $J_{\ell^*}$), which is disjoint from all segments
$I_\ell, J_\ell$ and contains $\gg (\log^2 n)|A_t|$ vertices of
$N(A_t)$. Let $I_{t+1}=S_i$. Then we can perform a rotation of the
points in $A_t$ as we did in Case 1. Note that in the end  we will
get a set $A_{t+1}$ which has size $\gg (\log^2 n)|A_t|$. On the
other hand,  the set $B_t$ may shrink  by a factor of at most
$2(t+1)$. However, as we already mentioned in the beginning of the induction step, 
this loss in the size of $B_t$  will be compensated later when
we rotate right end points.

\vspace{0.15cm} {\bf Case 3.}\, If cases 1 and 2 do not
occur, then there are less than $(\log^4 n)|A_t|/45$
vertices outside the set $\cup_{\ell}(I_{\ell}\cup J_{\ell})$ and also less than $(\log^4 n)|A_t|/45$ vertices 
inside the set $\cup_{\ell} J_{\ell}$. Since the size of $N(A_t)$ is at least
$(\log^4 n)|A_t|/15$, we have at least $(\log^4 n)|A_t|/45$ neighbors of $A_t$
inside $\cup_{\ell} I_{\ell}$. Hence there is an index $\ell^*$
such that $I_{\ell^*}$ contains at least $\Omega\big((\log^4
n)|A_t|/t\big)>(\log^3 n)|A_t|$ vertices from $N(A_t)$. Let $L$ be the
set of all indices $1 \leq \ell \leq t$ such that
$I_{\ell}=I_{\ell^*}$. Since for every $a \in A_t$ there are at
most $t$ broken edges $e_\ell(a)$ inside $I_{\ell^*}$, the total
number of broken edges inside this segment is at most $t|A_t|\ll
(\log n)|A_t| $. Therefore, there is a subset $D$ of at least
$(\log^3 n)|A_t|/2$ vertices from $N(A_t)$ inside $I_{\ell^*}$, such
that for every vertex $u \in D$, there is a segment of the path $P$
of length $\log n$ (constant length here is enough) with center in
$u$, which contains no broken edges. By definition, for every $u \in
D$ there is a vertex $a_u \in A_t$  adjacent to $u$. If there are
several such vertices, we fix one of them for each $u$. By
properties 1--4, the broken edges $\{e_\ell(a), \ell \in L\}$
appear in the same order in $I_{\ell^*}$ for all $a \in A_t$, and
partition it into  at most $t+1$ segments. So, there exists a
subset $D'\subset D$ of size at least $|D|/(t+1)\gg (\log^2 n)
|A_t|$, and an index $j$ such that every $u \in D'$ appears in the
$j$-th segment (from the left) of $I_{\ell^*}\setminus
\{e_\ell(a_u), \ell \in L\}$, where $a_u$ is the neighbor of $u$ in
$A_t$ which we fixed.

Since the order of edges $\{e_\ell(a)\}$ is the same for all $a$
we conclude that there are two indices $\ell_1$ and $\ell_2$ such
that every $u \in D'$ appears after edge $e_{\ell_1}(a_u)$ and
before edge $e_{\ell_2}(a_u)$. By properties 1--4,  there is an
index $1 \leq h \leq 2t+1$ such that $u$ is contained in the
segment $Q_h(a_u,b)$, and this holds for all $u \in D'$ and $b \in
B_t$. Without loss of generality,  suppose that $z_h=0$, i.e., all
paths $P(a,b)$ traverse the segment $Q_{h}(a,b)$ from left to
right (in the direction of path $P$). Let $a'$ be the left
neighbor of $u$ on the path $P$. Then we can rotate  $P(a,b)$,
using $u$ as a pivot and keeping $b$ fixed, to obtain a new path
$P(a',b)$. Put $a' \in A_{t+1}$, and note that for all $b \in B_t$,
we should  have the same new broken edge $(a',u)$ which we denote
$e_{t+1}(a')$. Also, for all $1 \leq \ell \leq t$, let
$e_{\ell}(a')=e_{\ell}(a)$. By our construction, $e_{t+1}(a')$ is
located after $e_{\ell_1}(a')$ and immediately before
$e_{\ell_2}(a')$. So the order in which the $2t+1$ edges
$\{e_{\ell}(a')\} \cup\{e_{\ell}(b)\}$ appear on the path $P$ do
not depend on the pair $(a',b)$. The rest the of analysis showing
that all new maximum paths $P(a',b)$ have the same permutation
$\pi'$ and vector $z'$  is exactly the same as in Case 1 and
omitted. Finally, set $I_{t+1}=I_{\ell^*}$ and note that
$|A_{t+1}|=|D'| \gg (\log^2 n) |A_t|$. This concludes the
description of the rotation procedure for the left end points.

As we already mentioned, the rotation of the right end points 
can be done similarly. Therefore, this also completes the proof of the induction step.
\hfill $\Box$

\subsection{Sparse pseudo-random graphs}

Theorem \ref{hamilton-cycle} can be also applied to show that sparse 
$(n,D,\lambda)$-graphs have local resilience with respect to hamiltonicity.
For this we first need to prove the following analog
of Lemma \ref{expansion}.

\begin{lemma}
\label{expansion-pseudo-random} For any fixed $\epsilon>0$ and
sufficiently large $n$ the following holds. Let $G$ be an $(n,D,
\lambda)$-graph with $D/\lambda > \log^2 n$ and let 
$H$ be a subgraph of $G$ with maximum degree $\Delta(H) \leq
(1/2-\epsilon)D$, then the graph $G'=G-H$ is

\begin{itemize}

\item connected;

\item  every subset $U$ of $G'$ of size at most $n/\log^4 n$
satisfies that $|N_{G'}(U)| \geq |U| \log^4 n/15$;

\item every subset $W$ of size at least $n/\log^3 n$ has
$|N_{G'}(W)| \geq \frac{1+\epsilon}{2}n$.

\end{itemize}
\end{lemma}

\noindent {\bf Proof.}\, We  use the following well known estimate
on the edge distribution of an $(n,D,\lambda)$-graph $G$ (see,
e.g., \cite{AloSpe00}, Corollary 9.2.5). For every two (not
necessarily disjoint) subsets $B,C\subseteq V$, let $e(B,C)$
denote the number of ordered pairs $(u,v)$ with $u \in B, v \in C$
such that $uv$ is an edge. Note that if $u,v \in B \cap C$, then
the edge $uv$ contributes $2$ to $e(B,C)$. In this notation,
$$
\left|e(B,C)-\frac{D}{n}|B||C|\right|\le \lambda \sqrt{|B||C|}\ .
$$

Let $U$ be a subset of $G'$ of size at most $n/\log^4 n$ and let
$X=U \cup N_{G'}(U)$. Observe  that since $\Delta(H) \leq
(1/2-\epsilon)D$ we have that the minimum degree of $G'$ is at
least $D-\Delta(H)>D/2$. Therefore
\begin{eqnarray*}
e_G(U,X) &\geq& e_{G'}(U,X) \geq \frac{D}{2}|U|-2e_G(U) \geq
\frac{D}{2}|U|-\frac{D}{n}|U|^2-\lambda|U|\\
& \geq& \left(\frac{D}{2}-\frac{D}{\log^4 n}-\frac{D}{\log^2 n}
\right)|U| \geq \frac{2}{5}D|U|.
\end{eqnarray*}
Suppose that $|X| \leq \frac{\log^4 n}{14}|U|\leq n/14$. Then
\begin{eqnarray*}
e_G(U,X) &\leq& \frac{D}{n}|U||X|+\lambda\sqrt{|U||X|} \leq
\frac{D}{14}|U|+\frac{\lambda \log^2 n}{\sqrt{14}}|U|\\
&<&
\left(\frac{D}{14}+\frac{D}{\sqrt{14}}\right)|U|<\frac{2}{5}D|U|.
\end{eqnarray*}
This contradiction implies that $|X|\geq \frac{\log^4 n}{14}|U|$ and
that $|N_{G'}(U)|\geq |X|-|U|>\frac{\log^4 n}{15}|U|$.

Now suppose that $W$ is a set of vertices of size at least
$n/\log^3 n$ such that $|N_{G'}(W)|<\frac{1+\epsilon}{2}n$. Let $Y
\subseteq W$ of size $|Y|=n/\log^3 n$. Then also
$|N_{G'}(Y)|<\frac{1+\epsilon}{2}n$ and there is a subset $Z$ of
$G'$ of size at least $\frac{1-\epsilon}{2}n$ such that
$e_{G'}(Y,Z)=0$. Hence all the edges of $G$ from $Y$ to $Z$ were
deleted when we removed $H$ and therefore $e_{G}(Y,Z)\leq
\Delta(H)|Y|\leq (1/2-\epsilon)D|Y|$. On the other hand
\begin{eqnarray*}
e_{G}(Y,Z) &\geq& \frac{D}{n}|Y||Z|-\lambda\sqrt{|Y||Z|} \geq
\frac{1-\epsilon}{2}D|Y|- \lambda \log^{3/2} n |Y|\\
&\geq& \left(\frac{1-\epsilon}{2}-\frac{1}{\sqrt{\log
n}}\right)D|Y|> \big(1/2-\epsilon\big)D|Y|.
\end{eqnarray*}
This contradiction implies that $|N_{G'}(W)|\geq \frac{1+\epsilon}{2}n$. In particular,
every connected component of $G'$ has size larger than $n/2$ so $G'$ is connected.
\hfill $\Box$

\vspace{0.25cm} This lemma and Theorem \ref{hamilton-cycle} imply
the following theorem, which strengthens the result of Krivelevich
and Sudakov \cite{KS} on the hamiltonicity of   sparse
$(n,D,\lambda)$-graphs.

\begin{theo}
\label{hamilton-pseudo-random} For any fixed $\epsilon>0$ and $n$
sufficiently large the following holds. Let $G$ be an
$(n,D,\lambda)$-graph such that $D/\lambda>\log^2 n$. If $H$ is a
subgraph of $G$ with maximum degree $\Delta(H) \leq
(1/2-\epsilon)D$ then $G'=G-H$ contains a Hamilton cycle. In other
words, the local resilience of $G$ with respect to hamiltonicity is
at least $(1/2-\ep)D$.
\end{theo}

\vspace{0.25cm} \noindent {\bf Remark.}\, Using a more careful
analysis one can show that in the statement of Proposition 
\ref{hamilton-cycle} it suffices to assume that  small
sets expand by a factor of $\log^{2+\delta} n$ for arbitrary fixed
$\delta >0$ (instead of $\log^4 n/15$). Therefore the assertion of
Theorem \ref{hamilton-pseudo-random} holds already when
$p \geq \log^{2+\delta} n/n$ and in
Theorem \ref{hamilton-pseudo-random} it is enough to assume that
$D/\lambda>\log^{1+\delta} n$.

\section{Chromatic number}
In this section we obtain results on local resilience  of random
graphs with respect to colorability. We start with the proof of
Theorem \ref{t1}, which states that for every positive integer $d$
and for every $\epsilon>0$ there exist a constant $c(d,\epsilon)$
such that if the edge probability $p>c/n$ then almost surely
$$\max_{H,\, \Delta(H)\leq d}~ \chi\big(G(n,p)\cup
H\big)\leq (1+\epsilon)\chi(G(n,p)).$$

Since for dense random graphs we obtain much stronger result, we
prove the theorem only for $p=o(1)$. To do so we will need the
following two lemmas which summarize some useful properties of
random graphs. The first one is well known and gives the
asymptotic value of the chromatic number of random graph. It is an
immediate corollary of the result of \L uczak (mentioned earlier)
together with a result of Shamir and Spencer about concentration
of the chromatic number of $G(n,p)$. We refer an interested reader
to the Chapter 7 of \cite{JanLucRuc00} for the detailed
description of these results.

\begin{lemma}
\label{le1}
For every $\epsilon>0$ there exist $c(\epsilon)$ such that if
$c/n<p=o(1)$ then with probability at least $1-n^{-1}$
$$(1-\epsilon/4)\frac{np}{2\log np}<\chi(G(n,p))<(1+\epsilon/4)\frac{np}{2\log np}.$$
\end{lemma}

A graph is {\em $d$-degenerate} if every subgraph of it
contains a vertex of degree at most $d$. It is a simple and
well known fact that every $d$-degenerate graph is $(d+1)$-colorable.
We will also need the following estimate of the density of small subgraphs
of $G(n,p)$. This estimate is standard (see, e.g., \cite{{JanLucRuc00}})
and can be easily verified using the first moment method.

\begin{lemma}
\label{le2}
If $ p \geq 200/n $, then a.s. every
$t \leq n/{\log}^{2} (np)$ vertices of random graph $G(n,p)$
span fewer than $(2np/{\log}^{2} (np))t$ edges.
Therefore any subgraph of this graph
induced by a subset $U$ of
size $|U|\leq n/{\log}^{2}(np)$, is $4np/{\log}^{2} (np)$-degenerate.
\end{lemma}

\noindent {\bf Proof of Theorem \ref{t1}.}\, Let $H$ be a graph on
the vertex set $[n]$ with $\Delta(H)\leq d$ and let $s=2d\log^2
(np)$. Pick a partition of vertices of $H$ into parts $V_1,
\ldots, V_s$ of size $n/s$ uniformly at random. Denote by $Y$ the
number of edges of $H$ which fall into one of the parts $V_i$. For
every edge, the probability of it falls into one of the $V_i$  is
$(1+o(1))/s$. By linearity of expectation  $\mathbb{E}[Y]\leq
(1+o(1))\frac{nd}{2s}$. By Markov's inequality, with probability
at least $(1+o(1))/2>1/3$ (i.e. for at least $1/3$ of all possible
partitions) we have that $Y \leq nd/s=\frac {n}{2\log^2 (np)}$.
Hence for every graph $H$, for at least $1/3$ of all possible
partitions, the number of edges of $H$ falling into one of the
parts is at most $\frac {n}{2\log^2 (np)}$.

Now fix a  partition $V_1, \ldots, V_s$ and consider the
probability of the event that for some $i$, the subgraph of
$G=G(n,p)$ induced by $V_i$ has chromatic number greater than
$(1+\epsilon/4)\frac{np/s}{2\log (np/s)}$. If $np$ is sufficiently
large, then by Lemma \ref{le1} the probability of this event is at
most $s\cdot (n/s)^{-1}=s^2/n=o(1)$. Therefore with probability
$1-o(1)>2/3$ (i.e. for almost every partition) we have that the
chromatic number of $G[V_i]$ is bounded by
$(1+\epsilon/4)\frac{np/s}{2\log (np/s)}$. This implies that
almost surely for random graph $G=G(n,p)$ and every graph $H$ with
maximum degree $d$ there is a partition $V_1, \ldots, V_s$ such
that there are at most $\frac {n}{2\log^2 (np)}$ edges of $H$
which are contained in some part $V_i$ and for every $i$ we have
that $\chi(G[V_i]) \leq (1+\epsilon/4)\frac{np/s}{2\log (np/s)}$.

Now we can color $G\cup H$ as follows. First color each induced
subgraph $G[V_i]$ by at most $(1+\epsilon/4)\frac{np/s}{2\log
(np/s)}$ colors, using new colors for every $i$. This may create
some monochromatic edges, all of which belong to $H$ and are
contained in one of the $V_i$. The number of monochromatic edges
after this process  is at most  $\frac{n}{2\log^2 (np)}$. These
edges are adjacent to a set of at most $n/ \log^2 (np)$ vertices
which we denote by $U$. By Lemma \ref{le2}, $U$ induces
$4np/{\log}^{2} (np)$-degenerate subgraph of $G=G(n,p)$. Since the
maximum degree of $H$ is at most $d$, the subgraph of $G\cup H$
induced by $U$ is $(4np/{\log}^{2} (np)+d)$-degenerate. Thus we
can color it by $4np/{\log}^{2} (np)+d+1$ fresh colors. This
implies that
\begin{eqnarray*}
\chi\big(G(n,p)\cup H\big)&\leq& s (1+\epsilon/4)\frac{np/s}{2\log (np/s)}
+\frac{4np}{\log^{2} (np)}+d+1\\
&=&
(1+\epsilon/4)\frac{np}{2\log (np)-4\log \log (np)-\log (2d)}+
\frac{4np}{\log^{2} (np)}+d+1.
\end{eqnarray*}
By choosing an appropriate constant $c(d,\epsilon)$ we see that for $p>c/n$
the right hand side of the last inequality can be made smaller than
$(1+\epsilon)\chi(G(n,p))$. \hfill $\Box$

It is easy to see from the proof that if
$p=n^{-\alpha}$ for some constant $\alpha<1$ then
the assertion of the Theorem \ref{t1} remains valid even
if we allow to the maximum degree of $H$ to be as large as $e^{\log^{1-\delta}n},
\delta>0$. On the other hand for very dense random graphs we can significantly
improve this theorem as follows (this is a quantitative version of Theorem 2.6).

\begin{theo}
\label{t2}
Let $0<\epsilon<1/3$ be a constant. If the edge probability $p(n)$ satisfies
$n^{-1/3+\epsilon} \leq p(n) \leq 3/4$ then almost surely
$$\max_{H} \chi\big(G(n,p)\cup
H\big)=(1+o(1))\chi(G(n,p))=(1+o(1))\frac{n}{2\log_b(np)},$$
where $b=1/(1-p)$ and the maximum is taken over all graphs with
$\Delta(H) \leq \frac{np^2}{\log^5 n}$.
\end{theo}

When $p$ is a constant this result implies that a.s. one will not
affect the asymptotic value of the chromatic number of $G(n,p)$
even by adding to it any graph with maximum degree $O(n/\log^5
n)$. The magnitude of  the maximal degree is clearly best possible
up to a polylogarithmic factor. We believe that a stronger result
holds for all values of edge probability and will state a
conjecture  in the concluding remarks.

To prove Theorem \ref{t2} we need to recall some additional properties of random
graphs.
Let $k_0=k_0(n,p)$ be defined by
$$k_0=\max\left\{k: {n\choose k}(1-p)^{{k\choose 2}}\geq
n^4\right\}\,.$$
One can show easily that $k_0$ satisfies $k_0\sim 2\log_b(np)$
with $b=1/(1-p)$. Also, it follows from known results on the
asymptotic value of the independence number of $G(n,p)$ (see,
e.g., \cite{JanLucRuc00}) that a.s. the difference between $k_0$
and the independence number of $G(n,p)$ is bounded by an absolute
constant, as long as $p(n)\gg n^{-1/2+\epsilon}$ for any positive
$\epsilon>0$.

Given a graph $G$ on $n$ vertices and an integer $k_0$, a
collection $\cal C$ of pairs of vertices of $G$ is called a {\em
cover} if every independent set of size $k_0$ in $G$ contains a
pair from $\cal C$. We set $X=X(G)$ to be the minimum size of a
cover in $G$. Note that, by definition, $X(G)$ is precisely the number of 
edges we need to add to graph $G$ in order to destroy all independent 
sets of size $k_0$. When $G$ is
distributed according to $G(n,p)$, the the minimum size of a
cover in $G$ becomes a random variable. The following properties of $X(G(n,p))$ were established
in \cite{KSVW}.

\begin{lemma}
\label{le3}
Let $0<\epsilon<1/3$ be a constant and let
$n^{-1/3+\epsilon} \leq p \leq 3/4$ then

\noindent
$(i)$\, $\mathbb{E}[X]=\Omega\left(\frac{n^2p^2}{\log^2 n}\right)$.

\noindent
$(ii)$\, For every $n^2p > t>0$, $Pr[X\le \mathbb{E}[X]-t]\leq e^{-\frac{t^2}{2n^2p}}$.
\end{lemma}

\noindent
{\bf Proof of Theorem \ref{t2}.}\, Let $G=G(n,p)$ and let
$H$ be graph with maximum degree $\Delta(H) \leq \frac{np^2}{\log^5 n}$.
First note that a.s. every subset $W$ of
vertices of $G\cup H$ of size at least $|W|\geq n/\log^2 n$ contains
an independent set of size $(1-o(1))2\log_b(np)$, where
$b=1/(1-p)$. Indeed, the number of edges of $H$ inside $W$ is bounded by
$$\Delta(H)|W|\leq  \frac{np^2}{\log^5 n}|W|\leq \frac{|W|^2p^2}{\log^3 n}
\ll \frac{|W|^2p^2}{\log^2 n}.$$ From Lemma \ref{le3} we obtain
that the probability that there exists $H$ that destroys all the
independent sets of size $(1-o(1))2\log_b(np)$ inside $W$ is less
than
$$2^n e^{-\Omega\left(\frac{(|W|^2p^2/\log^2 n)^2}{2n^2p}\right)}\leq
2^n e^{-n^{1+3\epsilon-o(1)}}=o(1).$$

Now we can color $G\cup H$ as follows. As long as the set of
uncolored vertices has size at least $n/\log^2 n$, we can pull out
an independent set of size $(1-o(1))2\log_b(np)$ and color it by a
new color. Let $V_0$ be the set of all uncolored vertices after
this process. $V_0$ has size at most $n/\log^2 n$ and therefore by
Lemma \ref{le2} it induces a $4np/{\log}^{2} (np)$-degenerate
subgraph of $G$. Since the maximum degree of $H$ is at most
$np^2/\log^5 n$, the subgraph of $G\cup H$ induced by $V_0$ is
$4np/{\log}^{2} (np)+\Delta(H)$-degenerate. Since $4np/{\log}^{2}
(np)+\Delta(H) \leq 5np/{\log}^{2} (np)$, we can color it by
$5np/{\log}^{2} (np) +1$ fresh colors. As $2\log_b(np)=O(\log
n/p)$, it follows that the total number of colors used in the
whole operation is at most
$$\hspace{2.5cm}
 (1+o(1))\frac{n}{2\log_b(np)}+\frac{5np}{{\log}^{2}
(np)} +1=(1+o(1))\frac{n}{2\log_b(np)} \hspace{2.5cm}  $$

\noindent concluding the proof.

\section{Concluding remarks} \label{section:remark}

\subsection{Local resilience  and universality}

Another popular model of random graphs (already mentioned a few
times in the paper) is the model of random regular graphs. Given
two parameters $d$ and $n$, we fix the vertex set $V= \{1, \dots,
n\}$ and consider the set of all simple $d$-regular graphs on $V$,
equipped with the uniform probability.

As mentioned in the introduction, the notion of local resilience
first came up in \cite{KV}. The goal of this paper was to
establish a universal theorem between the Erd\H os-R\'enyi random
graph $G(n,p)$ and random  regular graph with  approximately the
same density ($d \approx np$). It was showed, among others, that
if a property $\CP$ holds almost surely for $G(n,p)$ and has
unbounded local resilience, then it holds almost surely for a
random regular graph of approximately the same density. For more
details and precise statements, we refer to \cite{KV}. (In this
paper, the term ``tolerance" was used instead of ``resilience".)

It was conjectured that most ``natural" graph properties have
unbounded local resilience  given $p$ being sufficiently large.
The results in this paper show that for some important properties
this is indeed the case. They also shed some light on the question
whether the random regular graph model is monotone with respect to
$d$, see the last section of \cite{KV}.

\subsection{A few open problems}
We conclude this paper with two  open problems.

\begin{problem}
Prove that if $pn/\log n \rightarrow \infty.$ then a.s. the local
resilience  of $G(n,p)$ with respect to hamiltonicity is
$\big(1/2+o(1)\big)np$.
\end{problem}

\begin{problem} Is it true
that for all constant $\delta>0$ and $n^{-1+\delta}\leq p \leq 3/4$  random graph
$G(n,p)$ a.s. satisfies that
$$\max_{H} \chi\big(G(n,p)\cup H\big)=(1+o(1))\chi(G(n,p)),$$
where the maximum is taken over all graphs $H$ with maximum
degree $d=o(np/\log np)$.
\end{problem}

\vspace{0.2cm}
\noindent
{\bf Acknowledgment.}
A part of this work was carried out when both authors were visiting Microsoft Research at
Redmond, WA. We would like to thank the members of the Theory Group at
Microsoft Research for their hospitality and for creating a
stimulating research environment. We also like to thank A. Frieze, P. Loh, H. Nguyen, J. Vondr\'ak, 
and all three referees for crefully reading this manuscript and for very 
helpful remarks.

\end{document}